\documentclass{pnastwo}

\usepackage{amssymb,amsfonts,amsmath}
\usepackage{bigcircle}

\newcommand{\NN}{\mathbf{N}}
\newcommand{\PP}{\mathbf{P}}
\newcommand{\QQ}{\mathbf{Q}}
\newcommand{\RR}{\mathbf{R}}
\newcommand{\ZZ}{\mathbf{Z}}
\newcommand{\psl}[2]{\mathrm{PSL}_{#1}({#2})}
\newcommand{\se}{\subseteq}
\newcommand{\inv}{^{-1}}
\newcommand{\supp}{\mathrm{supp}}
\newcommand{\betti}{\beta_{(2)}}
\newcommand{\qed}{\hfill \ensuremath{\Box}}
\newtheorem{NM_lemma}[theorem]{Lemma}
\newtheorem{NM_remark}[theorem]{Remark}

\contributor{Submitted to Proceedings
of the National Academy of Sciences of the United States of America}
\url{www.pnas.org/cgi/doi/10.1073/pnas.0709640104}
\copyrightyear{2008}
\issuedate{Issue Date}
\volume{Volume}
\issuenumber{Issue Number}

\begin{document}



\title{Groups of piecewise projective homeomorphisms}





\author{Nicolas Monod\affil{1}{EPFL, 1015 Lausanne, Switzerland}}

\contributor{Submitted to Proceedings of the National Academy of Sciences
of the United States of America}

\maketitle

\begin{article}

\begin{abstract}
The group of piecewise projective homeomorphisms of the line provides straightforward torsion-free counter-examples to the so-called von Neumann conjecture. The examples are so simple that many additional properties can be established.
\end{abstract}






\section{Introduction}
In 1924, Banach and Tarski accomplished a rather paradoxical feat. They proved that a solid ball can be decomposed into five pieces which are then moved around and reassembled in such a way as to obtain \emph{two} balls identical with the original one~\cite{Banach-Tarski}. This wellnigh miraculous duplication was based on Hausdorff's 1914 work~\cite{Hausdorff14_article}. 

\smallskip
In his 1929 study of the Hausdorff--Banach--Tarski paradox, von Neumann introduced the concept of amenable groups~\cite{vonNeumann29}. Tarski readily proved that amenability is the \emph{only} obstruction to paradoxical decompositions~\cite{Tarski29,Tarski38}. However, the known paradoxes relied more prosaically on the existence of non-abelian free subgroups. Therefore, the main open problem in the subject remained for half a century to find non-amenable groups without free subgroups. Von Neumann's name was apparently attached to it by Day in the 1950s. The problem was finally solved around 1980: Ol$'$shanski\u{\i} proved the non-amenability of the Tarski monsters that he had constructed~\cite{Olshanskii79T,Olshanskii80T,Olshanskii80}; Adyan showed that his work on Burnside groups yields non-amenability~\cite{Adyan1979, Adyan83}. Finitely presented examples were constructed another twenty years later by Ol$'$shanski\u{\i}--Sapir~\cite{Olshanskii-Sapir}. There are several more recent counter-examples~\cite{Ershov08,Osin09,Schlage-Puchta}.

\bigskip
Given any subring $A<\RR$, we shall define a group $G(A)$ and a subgroup $H(A)<G(A)$ of piecewise projective transformations. Those will provide concrete, uncomplicated new examples with many additional properties. Perhaps ironically, our short proof of non-amenability ultimately relies on basic free groups of matrices, as in Hausdorff's 1914 paradox, even though the Tits alternative~\cite{Tits72} shows that the examples cannot be linear themselves.

\subsection*{Construction}
\begin{flushright}
\begin{minipage}[t]{0.75\linewidth}
\small
\itshape
I saw the pale student of unhallowed arts kneeling beside the thing he had put together.
\upshape
\begin{flushright}
Mary~Shelley, \emph{Frankenstein}\\
(introduction to the 1831 edition)
\end{flushright}
\end{minipage}
\end{flushright}

\smallskip
Consider the natural action of the group $\psl2{\RR}$ on the projective line $\PP^1=\PP^1(\RR)$. We endow $\PP^1$ with its $\RR$-topology making it a topological circle. We denote by $G$ the group of all homeomorphisms of $\PP^1$ which are piecewise in $\psl2{\RR}$, each piece being an interval of $\PP^1$, with finitely many pieces. We let $H<G$ be the subgroup fixing the point $\infty\in\PP^1$ corresponding to the first basis vector of $\RR^2$. Thus $H$ is left-orderable since it acts faithfully on the topological line $\PP^1\setminus\{\infty\}$, preserving orientations. It follows in particular that $H$ is torsion-free.

\medskip
Given a subring $A<\RR$, we denote by $P_A\se\PP^1$ the collection of all fixed points of all hyperbolic elements of $\psl2{A}$. This set is $\psl2{A}$-invariant and is countable if $A$ is so. We define $G(A)$ to be the subgroup of $G$ given by all elements that are piecewise in $\psl2{A}$ with all interval endpoints in $P_A$. We write $H(A)=G(A)\cap H$, which is the stabilizer of $\infty$ in $G(A)$.

\medskip
The main result of this article is the following, which relies on a new method for proving non-amenability.

\begin{theorem}\label{thm:nam}
The group $H(A)$ is non-amenable if $A\neq \ZZ$.
\end{theorem}

The next result is a sequacious generalization of the corresponding theorem of Brin--Squier about piecewise affine transformations~\cite{Brin-Squier} and we claim no originality.

\begin{theorem}\label{thm:no-free}
The group $H$ does not contain any non-abelian free subgroup. Thus, $H(A)$ inherits this property for any subring $A<\RR$.
\end{theorem}

Thus already $H=H(\RR)$ itself is a counter-example to the von Neumann conjecture. Writing $H(A)$ as the directed union of its finitely generated subgroups, we deduce:

\begin{corollary}
For $A\neq \ZZ$, the groups $H(A)$ contain finitely generated subgroups that are simultaneously non-amenable and without non-abelian free subgroups.
\end{corollary}

\subsection*{Further properties }
The groups $H(A)$ seem to enjoy a number of additional interesting properties, some of which are weaker forms of amenability. In the last section, we shall prove the following five propositions (and recall the terminology). Here $A<\RR$ is an arbitrary subring.

\begin{proposition}\label{prop:betti}
All $L^2$-Betti numbers of $H(A)$ and of $G(A)$ vanish.
\end{proposition}

\begin{proposition}\label{prop:inner}
The group $H(A)$ is inner amenable.
\end{proposition}

\begin{proposition}\label{prop:bio}
The group $H$ is bi-orderable and hence so are all its subgroups. It follows that there is no non-trivial homomorphism from any Kazhdan group to $H$.
\end{proposition}

\begin{proposition}\label{prop:coamen}
Let $E\se\PP^1$ be any subset. Then the subgroup of $H(A)$ which fixes $E$ pointwise is co-amenable in $H(A)$ unless $E$ is dense (in which case the subgroup is trivial).
\end{proposition}

\begin{proposition}\label{prop:cat0}
If $H(A)$ acts by isometries on any proper CAT(0) space, then either it fixes a point at infinity or it preserves a Euclidean subspace.
\end{proposition}

One can also check that $H(A)$ satisfies no group law and has vanishing properties in bounded cohomology (see below).

\section{Non-amenability}
An obvious difference between the actions of $\psl2A $ and of $H(A)$ on $\PP^1$ is that the latter group fixes $\infty$ whilst the former does not. The next proposition shows that this is the only difference as far as the orbit structure is concerned.

\begin{proposition}\label{prop:orbit}
Let $A<\RR$ be any subring and let $p\in \PP^1\setminus\{\infty\}$. Then
$$\psl2A \cdot p \ \se\ \{\infty\} \cup H(A)\cdot p.$$
Thus, the equivalence relations induced by the actions of $\psl2A $ and of $H(A)$ on $\PP^1$ coincide when restricted to $\PP^1\setminus \{\infty\}$.
\end{proposition}

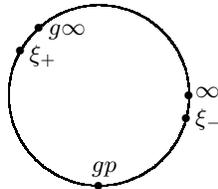
\begin{figure}[h]
\setlength{\unitlength}{0.3cm}
\begin{picture}(8,1)
\end{picture}
\begin{picture}(10,7.1)
\bigcircle{5.0}{4.0}{4.0}

\put(9,4){\circle*{0.3}}
\put(9.3,3.9){$\infty$}

\put(5,0){\circle*{0.3}}
\put(4.7,0.6){$g p$}

\put(2.354,7){\circle*{0.3}}
\put(2.8,6.7){$g \infty$}

\put(1.5359,6){\circle*{0.3}}
\put(1.95,5.5){${\xi_+}$}

\put(8.873,3){\circle*{0.3}}
\put(9.3,2.7){${\xi_-}$}
\end{picture}
\caption{The desired configuration of $\xi_\pm$}\label{fig:desired}
\end{figure}

\medskip
\emph{Proof.}
We need to show that given $g\in \psl2A$ with $g p\neq \infty$, there is an element $h\in H(A)$ such that $h p = g p$. We assume $g\infty \neq \infty$ since otherwise $h=g$ will do. Equivalently, we need an element $q\in G(A)$ fixing $g p$ and such that $q \infty = g\infty$, writing $h=q\inv g$. It suffices to find a hyperbolic element $q_0\in\psl2A$ with $q_0 \infty = g\infty$ and whose fixed points $\xi_\pm\in\PP^1$ separate $gp$ from both $\infty$ and $g \infty$, see Figure~\ref{fig:desired}. Indeed, we can then define $q$ to be the identity on the component of $\PP^1\setminus\{\xi_\pm\}$ containing $g p$, and define $q$ to coincide with $q_0$ on the other component.

Let $\begin{pmatrix}a&b\\ c&d\end{pmatrix}$ be a matrix representative of $g$; thus, $a,b,c,d\in A$ and $ad-bc=1$.  The assumption $g\infty \neq \infty$ implies $c\neq 0$ and thus we can assume $c>0$. Let $q_0$ be given by $\begin{pmatrix}a&b + r a\\ c&d + r c\end{pmatrix}$ with $r\in A$ to be determined later; thus $q_0\infty = g\infty$. This matrix is hyperbolic as soon as $|r|$ is large enough to ensure that the trace $\tau=a+d+r c$ is larger than $2$ in absolute value. We only need to show that a suitable choice of $r$ will ensure the above condition on $\xi_\pm$. Notice that $\infty$ and $g\infty$ lie in the same component of $\PP^1\setminus\{\xi_\pm\}$ since $q_0$ preserves these components and sends $\infty$ to $g\infty$. In conclusion, it suffices to prove the following two claims: (1)~as $|r|\to\infty$, the set $\{\xi_\pm\}$ converges to $\{\infty, g \infty\}$; (2)~changing the sign of $r$ (when $|r|$ is large) will change the component of $\PP^1\setminus\{\infty, g \infty\}$ in which $\xi_\pm$ lie (we need it to be the component of $g p$). The claims can be proved by elementary dynamical considerations; we shall instead verify them explicitly.

\smallskip
The fixed points $\xi_\pm$ are represented by the eigenvectors $\begin{pmatrix}x_\pm\\ c\end{pmatrix}$, where $x_\pm=\lambda_\pm-d-r c$ and where $\lambda_\pm= (\tau\pm\sqrt{\tau^2-4})/2$ are the eigenvalues. Now $\lim_{r\to+\infty} \lambda_+ = +\infty$ implies $\lim_{r\to+\infty} \lambda_- = 0$ since $\lambda_+\lambda_-=1$ and therefore $\lim_{r\to+\infty} x_-=-\infty$. Similarly, $\lim_{r\to-\infty} x_+ = +\infty$ (Figure~\ref{fig:desired} depicts the case $r>0$). This already proves claim~(2) and half of claim~(1). Since $g\infty=[a\colon\! c]$, it only remains to verify that both $\lim_{r\to+\infty} x_+$ and $\lim_{r\to-\infty} x_-$ converge to $a$, which is a direct computation.
\medskip
\qed

We recall that a measurable equivalence relation with countable classes is \emph{amenable} if there is an a.e.\ defined measurable assignment of a mean on the orbit of each point in such a way that the means of two equivalent points coincide. We refer e.g.\ to~\cite{Connes-Feldman-Weiss} and~\cite{Kechris-Miller04} for background on amenable equivalence relations. It follows from this definition that any relation produced by a measurable action of a (countable) amenable group is amenable, by push-forward of the mean~\cite[1.6(1)]{Schmidt90}. An a.e.\ free action of a countable group is amenable in Zimmer's sense~\cite[4.3]{Zimmer84} if and only if the associated relation is amenable; see~\cite[Thm.~A]{Adams-Elliott-Giordano}.

\medskip
\emph{Proof of Theorem~\ref{thm:nam}.}
Let $A\neq \ZZ$ be a subring of $\RR$. Then $A$ contains a countable subring $A'<A$ which is dense in $\RR$. Since $H(A')$ is a subgroup of $H(A)$, we can assume that $A$ itself is countable dense. Now $H(A)$ is a countable group and $\Gamma:=\psl2A$ is a countable dense subgroup of $\psl2{\RR}$.

It is proved in Th\'eor\`eme~3 of~\cite{Carriere-Ghys} that the equivalence relation on $\psl2{\RR}$ induced by the multiplication action of $\Gamma$ is non-amenable; see also Remarks~\ref{rem:pingpong} and~\ref{rem:alt} below. Equivalently, the $\Gamma$-action on $\psl2{\RR}$ is non-amenable. Viewing $\PP^1$ as a homogeneous space of $\psl2{\RR}$, it follows that the $\Gamma$-action on $\PP^1$ is non-amenable. Indeed, amenability is preserved under extensions, see~\cite[2.4]{Zimmer78b} or~\cite[Cor.~C]{Adams-Elliott-Giordano}. This action is a.e.\ free since any non-trivial element has at most two fixed points. Thus the relation induced by $\Gamma$ on $\PP^1$ is non-amenable. Restricting to $\PP^1\setminus \{\infty\}$, we deduce from Proposition~\ref{prop:orbit} that the relation induced by the $H(A)$-action is also non-amenable. (Amenability is preserved under restriction~\cite[9.3]{Kechris-Miller04}, but here $\{\infty\}$ is a null-set anyway.) Thus $H(A)$ is a non-amenable group.
\medskip
\qed

\begin{NM_remark}\label{rem:pingpong}
We recall from~\cite{Carriere-Ghys} that the non-amenability of the $\Gamma$-relation on $\psl2{\RR}$ is a general consequence of the existence of a non-discrete non-abelian free subgroup of $\Gamma$. Thus the main point of our appeal to~\cite{Carriere-Ghys} is the existence of this non-discrete free subgroup, but this is much easier to prove directly in the present case of $\Gamma=\psl2A$ than for general non-discrete non-soluble $\Gamma$.
\end{NM_remark}

\begin{NM_remark}\label{rem:alt}
Here is a direct argument avoiding all the above references in the examples of $A=\ZZ[\sqrt2]$ or $A=\ZZ[1/\ell]$, where $\ell$ is prime. We show directly that the $\Gamma$-action on $\PP^1$ is not amenable. We consider $\Gamma$ as a lattice in $L:=\psl2{\RR}\times\psl2{\RR}$ in the first case and in $L:=\psl2{\RR}\times\psl2{\QQ_\ell}$ in the second case, both times in such a way that the $\Gamma$-action on $\PP^1$ extends to the $L$-action factoring through the first factor. If the $\Gamma$-action on $\PP^1$ were amenable, so would be the $L$-action (by co-amenability of the lattice). But of course $L$ does not act amenably since the stabilizer of any point contains the (non-amenable) second factor of $L$.
\end{NM_remark}

The non-discreteness of $A$ was essential in our proof, thus excluding $A=\ZZ$.

\begin{problem}
Is $H(\ZZ)$ amenable?
\end{problem}

The group $H(\ZZ)$ is related to Thompson's group $F$, for which the question of (non-)amenability is a notorious open problem. Indeed $F$ seems to be historically the first candidate for a counter-example to the so-called von Neumann conjecture. The relation is as follows: if we modify the definition of $H(\ZZ)$ by requiring that the breakpoints be rational, then all its elements are automatically $C^1$ and the resulting group is conjugated to $F$. The corresponding relation holds between $G(\ZZ)$ and Thompson's group $T$. These facts are attributed to a remark of Thurston around 1975 and a very detailed exposition can be found in~\cite{Martin_PhD}.

\section{$H$ is a free group free group}
We shall largely follow~\cite[\S\,3]{Brin-Squier}, the main difference being that we replace commutators by a non-trivial word in the \emph{second} derived subgroup of a free group on two generators.

\medskip
The \emph{support} $\supp(g)$ of an element $g\in H$ denotes the set $\{p: g p \neq p\}$, which is a finite union of open intervals. Any subgroup of $H$ fixing some point $p\in\PP^1$ has two canonical homomorphisms to the metabelian stabilizer of $p$ in $\psl2{\RR}$ given by left and right germs. Therefore, we deduce the following elementary fact, wherein $\langle f, g\rangle$ denotes the subgroup of $H$ generated by $f$ and $g$.

\begin{NM_lemma}\label{lem:germs}
If $f,g\in H$ have a common fixed point $p\in\PP^1$, then any element of the second derived subgroup $\langle f, g\rangle''$ acts trivially on a neighbourhood of $p$.\qed
\end{NM_lemma}

Theorem~\ref{thm:no-free} is an immediate consequence of the following more precise statement.

\begin{theorem}
Let $f,g\in H$. Either $\langle f, g\rangle$ is metabelian or it contains a free abelian group of rank two.
\end{theorem}

\medskip
\emph{Proof.}
We suppose that $\langle f, g\rangle$ is not metabelian, so that there is a word $w$ in the second derived subgroup of a free group on two generators such that $w(f,g)\in H$ is non-trivial. We now follow faithfully the proof of Theorem~3.2 in~\cite{Brin-Squier}, replacing $[f,g]$ by $w(f,g)$. For the reader's convenience, we sketch the argument; the details are on page~495 of~\cite{Brin-Squier} (or ~\cite[p.~232]{Cannon-Floyd-Parry}). Applying Lemma~\ref{lem:germs} to all endpoints $p$ of the connected components of $\supp(f)\cup\supp(g)$, we deduce that the closure of $\supp(w(f,g))$ is contained in $\supp(f)\cup\supp(g)$. This implies that some element of $\langle f, g\rangle$ will send any connected component of $\supp(w(f,g))$ to a disjoint interval. The needed element might depend on the connected component. However, upon replacing $w(f,g)$ by another non-trivial element $w_1\in \langle f, g\rangle''$ with minimal number of intersecting components with $\supp(f)\cup\supp(g)$, some element $h$ of $\langle f, g\rangle$ sends the whole of $\supp(w_1)$ to a set disjoint from it. The corresponding conjugate $w_2:=h w_1 h\inv$ will commute with $w_1$ and indeed these two elements generate freely a free abelian group.
\medskip
\qed

As pointed out to us by Cornulier, the above argument can be pushed so that $w_1$ and $h$ generate a wreath product $\ZZ\wr\ZZ$, compare~\cite[Thm.~21]{Guba-Sapir99} for the piecewise linear case.

\section{Lagniappe}
\emph{Proof of Proposition~\ref{prop:betti}.}
We refer to~\cite{Cheeger-Gromov} for the $L^2$-Betti numbers $\betti^n$, $n\in\NN$. Fix a large integer $n$ and let $\Gamma=G(H)$ or $H(A)$. Choose a set $F\se P_A$ of $n+1$ distinct points and let $\Lambda<\Gamma$ be the pointwise stabilizer of $F$. Any intersection $\Lambda^*$ of any (finite) number of conjugates of $\Lambda$ is still the pointwise stabilizer of a finite set $F^*$ containing $m\geq n+1$ points. The definition of $G(A)$ shows that $\Lambda^*$ is the product of $m$ infinite groups. The K\"unneth formula~\cite[\S\,2]{Cheeger-Gromov} implies $\betti^i(\Lambda^*)=0$ for all $i=0, \ldots, m-1$. In this situation, Theorem~1.3 of~\cite{Bader-Furman-Sauer} asserts $\betti^i(\Gamma)=0$ for all $i\leq m-1$.
\medskip
\qed

A subgroup $K$ of a group $J$ is called \emph{co-amenable} if there is an $J$-invariant mean on $J/K$. Equivalent characterizations, generalizations and unexpected examples can be found in~\cite{Eymard72} and~\cite{Monod-Popa}.

\smallskip

Recall that a group $J$ is \emph{inner amenable} if there is a conjugacy-invariant mean on $J\setminus\{e\}$. It is equivalent to exhibit such a mean that is invariant under the second derived subgroup $J''$ since the latter is co-amenable in $J$. Thus, Proposition~\ref{prop:inner} is a consequence of the stronger fact that $H(A)$ is ``\{asymptotically commutative\}-by-metabelian'' in a sense inspired by~\cite{Zeller-Meier} as follows.

\begin{proposition}\label{prop:inner:bis}
Let $A<\RR$ be any subring. For any finite set $S\se H(A)''$ there is a non-trivial element $h_S\in H(A)$ commuting with each element of $S$.
\end{proposition}

Indeed, any accumulation point of this net of point-masses at $h_S$ is $H(A)''$-invariant.

\medskip
\emph{Proof of Proposition~\ref{prop:inner:bis}.}
By the argument of Lemma~\ref{lem:germs}, there is a neighbourhood of $\infty$ on which all elements of $S$ are trivial. Thus is suffices to exhibit a non-trivial element $h_S$ of $H(A)$ which is supported in this neighbourhood. Notice that $\psl2{\ZZ}$ contains hyperbolic elements with both fixed points $\xi_\pm$ arbitrarily close to $\infty$, and on the same side. For instance, conjugate $\begin{pmatrix}2&1\\ 1&1\end{pmatrix}$ by $\begin{pmatrix}1&n\\0&1\end{pmatrix}$ for sufficiently large $n\in\NN$. We choose such an element $h_0$ with $\xi_\pm$ in the given neighbourhood and define $h_S$ to be trivial on the component of $\PP^1\setminus\{\xi_\pm\}$ containing $\infty$ and to coincide with $h_0$ on the other component.
\medskip
\qed

A group is called \emph{bi-orderable} if it carries a bi-invariant total order. The construction below is completely standard, compare e.g.~\cite[p.~233]{Cannon-Floyd-Parry} for a first-order version of our second-order argument.

\medskip
\emph{Proof of Proposition~\ref{prop:bio}.}
Choose an orientation of $\PP^1\setminus\{\infty\}$ and define a (right) germ at a point $p$ to be positive if either its first derivative is~$>1$ or if it is~$=1$ but the second derivative is~$>0$. Then define the set $H_+$ of positive elements of $H$ to consist of all transformations whose first non-trivial germ (starting from $\infty$ along the orientation) is positive. Now $H_+$ is a conjugacy invariant sub-semigroup and $H\setminus\{e\}$ is $H_+\sqcup H_+\inv$; this means that $H_+$ defines a bi-invariant total order.

Suppose now that we are given a homomorphism from a Kazhdan group to $H$. Its image is then a Kazhdan subgroup $K<H$. Kazhdan's property implies that $K$ is finitely generated. It has been known for a long time that any non-trivial finitely generated bi-orderable group has a non-trivial homomorphism to $\RR$: this follows ultimately from H\"older's 1901 work~\cite{Holder01} by looking at maximal convex subgroups and is explained in~\cite[\S\,2]{Kopytov-Medvedev}. But this is impossible for a Kazhdan group.
\medskip
\qed

\begin{NM_lemma}\label{lem:contr}
For any $p\in\PP^1\setminus\{\infty\}$ there is a sequence $\{g_n\}$ in $H(\ZZ)$ such that $g_n q$ converges to $\infty$ uniformly for $q$ in compact subsets of $\PP^1\setminus\{p\}$.
\end{NM_lemma}

\medskip
\emph{Proof.}
It suffices to show that for any open neighbourhoods $U$ and $V$ of $p$ and $\infty$ respectively in $\PP^1$, there is $g\in H(\ZZ)$ which maps $\PP^1\setminus U$ into $V$. Since the collection of pairs of fixed points of hyperbolic elements of $\psl2{\ZZ}$ is dense in $\PP^1\times \PP^1$, we can find hyperbolic matrices $h_1, h_2\in \psl2{\ZZ}$ with repelling fixed points $r_i$ in $U\setminus \{p\}$ and attracting fixed points $a_i$ in $V\setminus \{\infty\}$ and such that the cyclic order is $\infty, a_1, r_1, p, r_2, a_2$. Now we define $g$ to be a sufficiently high power of $h_1$ on the interval $[a_1, r_1]$ (for the above cyclic order), of $h_2$ on the interval $[r_2, a_2]$ and the identity elsewhere.
\medskip
\qed

\medskip
\emph{Proof of Proposition~\ref{prop:coamen}.}
Let $K$ be the pointwise stabilizer of a non-dense subset $E\se \PP^1$; it suffices to find a mean invariant under $H(A)''$. Let $\{g_n\}$ be the sequence provided by Lemma~\ref{lem:contr} for $p$ an interior point of the complement of $E$. Any accumulation point of the sequence of point-masses at $g_n K$ in $H(A)/K$ will do. Indeed, since any $g\in H(A)''$ is trivial in a neighbourhood of $\infty$, we have $g_n\inv g g_n\in K$ for $n$ large enough.
\medskip
\qed

The existence of two (or more) \emph{commuting} co-amenable subgroups is also a weak form of amenability. It is the key in the argument cited below.

\medskip
\emph{Proof of Proposition~\ref{prop:cat0}.}
Consider two disjoint non-empty open sets in $\PP^1$. The pointwise stabilizers of their complement commute with each other and are co-amenable by Proposition~\ref{prop:coamen}. In this situation, Corollary~2.2 of~\cite{Caprace-Monod_discrete} yields the desired conclusion.
\medskip
\qed

The properties used in this section show immediately that $H(A)$ fulfills the criterion of~\cite[Thm.~1.1]{Abert_laws} and thus satisfies no group law.

\medskip
Combining Theorems~\ref{thm:nam} and~\ref{thm:no-free} with the main result of~\cite{Monod-Ozawa2009}, we conclude that the wreath product $\ZZ\wr H$ is a torsion-free non-unitarisable group without free subgroups. We can replace it by a finitely generated subgroup upon choosing a non-amenable finitely generated subgroup of $H$. This provides some new examples towards Dixmier's problem, unsolved since 1950~\cite{Day50,Dixmier50,Nakamura-Takeda51}.

\medskip
Finally, we mention that our argument from Proposition~6.4 in~\cite{MonodVT} applies to show that the bounded cohomology $\mathrm{H}_\mathrm{b}^n(H(A), V)$ vanishes for all $n\in\NN$ and all mixing unitary representations $V$. More generally, it applies to any semi-separable coefficient module $V$ unless all finitely generated subgroups of $H(A)''$ have invariant vectors in $V$ (see~\cite{MonodVT} for details and definitions). This should be contrasted with the fact that amenability is characterized by the vanishing of bounded cohomology with all dual coefficients.





\begin{acknowledgments}
Supported in part by the ERC and the Swiss National Science Foundation.
\end{acknowledgments}


\begin{thebibliography}{10}

\bibitem{Banach-Tarski}
S.~Banach and A.~Tarski.
\newblock {Sur la d{\'e}composition des ensembles de points en parties
  respectivement congruentes [On decomposing point-sets into congruent parts].}
\newblock {\em Fund. math.}, 6:244--277, 1924.
\newblock French.

\bibitem{Hausdorff14_article}
F.~Hausdorff.
\newblock Bemerkung {\"u}ber den {I}nhalt von {P}unktmengen [Remarks on the content of point-sets].
\newblock {\em Math. Ann.}, 75(3):428--433, 1914.
\newblock German.

\bibitem{vonNeumann29}
J.~von Neumann.
\newblock Zur allgemeinen {T}heorie des {M}a{\ss}es [On the general theory of measure].
\newblock {\em Fund. Math.}, 13:73--116, 1929.
\newblock German.

\bibitem{Tarski29}
A.~Tarski.
\newblock Sur les fonctions additives dans les classes abstraites et leur
  application au probl\`eme de la mesure [On additive functions on abstract classes and applications to the measure problem].
\newblock {\em C. R. Soc. Sc. Varsovie}, 22:114--117, 1929.
\newblock French.

\bibitem{Tarski38}
A.~Tarski.
\newblock Algebraische {F}assung des {M}a{\ss}problems [An Algebraic formulation of the measure problem].
\newblock {\em Fund. Math.}, 31:47--66, 1938.
\newblock German.

\bibitem{Olshanskii79T}
A.~Y. Ol{$'$}shanski\u{\i}.
\newblock An infinite simple torsion-free {N}oetherian group.
\newblock {\em Izv. Akad. Nauk SSSR Ser. Mat.}, 43(6):1328--1393, 1979.

\bibitem{Olshanskii80T}
A.~Y. Ol{$'$}shanski\u{\i}.
\newblock An infinite group with subgroups of prime orders.
\newblock {\em Izv. Akad. Nauk SSSR Ser. Mat.}, 44(2):309--321, 479, 1980.

\bibitem{Olshanskii80}
A.~Y. Ol{$'$}shanski\u{\i}.
\newblock On the question of the existence of an invariant mean on a group.
\newblock {\em Uspekhi Mat. Nauk}, 35(4(214)):199--200, 1980.

\bibitem{Adyan1979}
S.~I. Adyan.
\newblock {\em The {B}urnside problem and identities in groups}, volume~95 of
  {\em Ergebnisse der Mathematik und ihrer Grenzgebiete}.
\newblock Springer-Verlag, Berlin, 1979.

\bibitem{Adyan83}
S.~I. Adyan.
\newblock Random walks on free periodic groups.
\newblock {\em Izv. Akad. Nauk SSSR Ser. Mat.}, 46(6):1139--1149, 1343, 1982.

\bibitem{Olshanskii-Sapir}
A.~Y. Ol{$'$}shanskii and M.~V. Sapir.
\newblock Non-amenable finitely presented torsion-by-cyclic groups.
\newblock {\em Publ. Math. Inst. Hautes \'Etudes Sci.}, 96:43--169 (2003),
  2002.

\bibitem{Ershov08}
M.~Ershov.
\newblock Golod-{S}hafarevich groups with property {$(T)$} and {K}ac-{M}oody
  groups.
\newblock {\em Duke Math. J.}, 145(2):309--339, 2008.

\bibitem{Osin09}
D.~V. Osin.
\newblock {$L^2$}-{B}etti numbers and non-unitarizable groups without free
  subgroups.
\newblock {\em Int. Math. Res. Not. IMRN}, (22):4220--4231, 2009.

\bibitem{Schlage-Puchta}
J.-C. Schlage-Puchta.
\newblock A {$p$}-group with positive rank gradient.
\newblock {\em J. Group Theory}, 15(2):261--270, 2012.

\bibitem{Tits72}
J.~Tits.
\newblock Free subgroups in linear groups.
\newblock {\em J. Algebra}, 20:250--270, 1972.

\bibitem{Brin-Squier}
M.~G. Brin and C.~C. Squier.
\newblock Groups of piecewise linear homeomorphisms of the real line.
\newblock {\em Invent. Math.}, 79(3):485--498, 1985.

\bibitem{Connes-Feldman-Weiss}
A.~Connes, J.~Feldman, and B.~Weiss.
\newblock An amenable equivalence relation is generated by a single
  transformation.
\newblock {\em Ergodic Theory Dynamical Systems}, 1(4):431--450 (1982), 1981.

\bibitem{Kechris-Miller04}
A.~S. Kechris and B.~D. Miller.
\newblock {\em Topics in orbit equivalence}, volume 1852 of {\em Lecture Notes
  in Mathematics}.
\newblock Springer-Verlag, Berlin, 2004.

\bibitem{Schmidt90}
K.~Schmidt.
\newblock {\em Algebraic ideas in ergodic theory}, volume~76 of {\em {CBMS}
  Regional Conference Series in Mathematics}.
\newblock Published for the Conference Board of the Mathematical Sciences, Washington, DC, 1990.

\bibitem{Zimmer84}
R.~J. Zimmer.
\newblock {\em Ergodic theory and semisimple groups}.
\newblock Birkh\"auser Verlag, Basel, 1984.

\bibitem{Adams-Elliott-Giordano}
S.~Adams, G.~A. Elliott, and T.~Giordano.
\newblock Amenable actions of groups.
\newblock {\em Trans. Amer. Math. Soc.}, 344(2):803--822, 1994.

\bibitem{Carriere-Ghys}
Y.~Carri{\`e}re and {\'E}.~Ghys.
\newblock Relations d'\'equivalence moyennables sur les groupes de {L}ie [Amenable equivalence relations on Lie groups].
\newblock {\em C. R. Acad. Sci. Paris S\'er. I Math.}, 300(19):677--680, 1985.
\newblock French.

\bibitem{Zimmer78b}
R.~J. Zimmer.
\newblock Amenable ergodic group actions and an application to {P}oisson
  boundaries of random walks.
\newblock {\em J. Functional Analysis}, 27(3):350--372, 1978.

\bibitem{Martin_PhD}
X.~Martin.
\newblock {\em Sur la g{\'e}om{\'e}trie du groupe de {T}hompson [On the geometry of {T}hompson's group].}
\newblock PhD thesis, Universit{\'e} de {G}renoble, 2002.
\newblock French.

\bibitem{Cannon-Floyd-Parry}
J.~W. Cannon, W.~J. Floyd, and W.~R. Parry.
\newblock Introductory notes on {R}ichard {T}hompson's groups.
\newblock {\em Enseign. Math. (2)}, 42(3-4):215--256, 1996.

\bibitem{Guba-Sapir99}
V.~S. Guba and M.~V. Sapir.
\newblock On subgroups of the {R}. {T}hompson group {$F$} and other diagram
  groups.
\newblock {\em Mat. Sb.}, 190(8):1077--1130, 1999.

\bibitem{Cheeger-Gromov}
J.~Cheeger and M.~Gromov.
\newblock {$L\sb 2$}-cohomology and group cohomology.
\newblock {\em Topology}, 25(2):189--215, 1986.

\bibitem{Bader-Furman-Sauer}
U.~Bader, A.~Furman, and R.~Sauer.
\newblock Weak notions of normality and vanishing up to rank in
  {$L^2$}-cohomology.
\newblock preprint {\tt arXiv:1206:4793v1}, 2012.

\bibitem{Eymard72}
P.~Eymard.
\newblock {\em Moyennes invariantes et repr\'esentations unitaires [Invariant means and unitary representations].}
\newblock Springer-Verlag, Berlin, 1972.
\newblock Lecture Notes in Mathematics, Vol. 300.
\newblock French.

\bibitem{Monod-Popa}
N.~Monod and S.~Popa.
\newblock On co-amenability for groups and von {N}eumann algebras.
\newblock {\em C. R. Math. Acad. Sci. Soc. R. Can.}, 25(3):82--87, 2003.

\bibitem{Zeller-Meier}
G.~Zeller-Meier.
\newblock Deux nouveaux facteurs de type {${\rm II}_{1}$} [Two new {${\rm II}_{1}$} factors].
\newblock {\em Invent. Math.}, 7:235--242, 1969.
\newblock French.

\bibitem{Holder01}
O.~H{\"o}lder.
\newblock {Die {A}xiome der {Q}uantit{\"a}t und die {L}ehre vom {M}a{\ss} [The axioms of quantity and the study of measure].}
\newblock {\em Leipz. Ber.}, 53:1--64, 1901.
\newblock German.

\bibitem{Kopytov-Medvedev}
V.~M. Kopytov and N.~Y. Medvedev.
\newblock {\em Right-ordered groups}.
\newblock Consultants Bureau, New York, 1996.

\bibitem{Caprace-Monod_discrete}
P.-E. Caprace and N.~Monod.
\newblock {Isometry groups of non-positively curved spaces: discrete
  subgroups}.
\newblock {\em J Topology}, 2(4):701--746, 2009.

\bibitem{Abert_laws}
M.~Ab{\'e}rt.
\newblock Group laws and free subgroups in topological groups.
\newblock {\em Bull. London Math. Soc.}, 37(4):525--534, 2005.

\bibitem{Monod-Ozawa2009}
N.~Monod and N.~Ozawa.
\newblock The {D}ixmier problem, lamplighters and {B}urnside groups.
\newblock {\em J. Funct. Anal.}, 1:255--259, 2010.

\bibitem{Day50}
M.~M. Day.
\newblock Means for the bounded functions and ergodicity of the bounded
  representations of semi-groups.
\newblock {\em Trans. Amer. Math. Soc.}, 69:276--291, 1950.

\bibitem{Dixmier50}
J.~Dixmier.
\newblock Les moyennes invariantes dans les semi-groupes et leurs applications [Invariant means on semi-groups and applications].
\newblock {\em Acta Sci. Math. Szeged}, 12:213--227, 1950.
\newblock French.

\bibitem{Nakamura-Takeda51}
M.~Nakamura and Z.~Takeda.
\newblock Group representation and {B}anach limit.
\newblock {\em T\^ohoku Math. J. (2)}, 3:132--135, 1951.

\bibitem{MonodVT}
N.~Monod.
\newblock On the bounded cohomology of semi-simple groups, {$S$}-arithmetic
  groups and products.
\newblock {\em J. Reine Angew. Math. [Crelle's J.]}, 640:167--202, 2010.

\end{thebibliography}




\end{article}








\end{document}